\documentclass[12pt]{amsart}
\usepackage{amssymb}
\usepackage{epsf}

\newtheorem{thm}{Theorem}
\newtheorem{prop}[thm]{Proposition}
\newtheorem{cor}[thm]{Corollary}

\begin{document}
\title[Real Rational Curves in Grassmannians]{Real
Rational Curves in Grassmannians}
\author{Frank Sottile}
\address{\hskip-\parindent
	Frank Sottile\\
	Department of Mathematics\\
        University of Wisconsin\\
        Van Vleck Hall\\
        480 Lincoln Drive\\
        Madison, Wisconsin 53706-1388\\
        USA}
\email{sottile@math.wisc.edu}
\urladdr{http://www.math.wisc.edu/\~{}sottile}
\date{23 April 1999}
\thanks{Research at MSRI supported in part by by NSF grant DMS-9701755}

\begin{abstract}
Fulton asked how many solutions to a problem of enumerative geometry can be
real, when that problem is one of counting geometric figures of some kind
having specified position with respect to some general fixed figures.  For
the problem of plane conics tangent to five general conics, the (surprising)
answer is that all 3264 may be real.  Similarly, given any problem of
enumerating $p$-planes incident on some general fixed subspaces, there are
real fixed subspaces such that each of the (finitely many) incident
$p$-planes are real.  We show that the problem of enumerating parameterized
rational curves in a Grassmannian satisfying simple (codimension 1)
conditions may have all of its solutions be real.
\end{abstract}

\maketitle

\section*{Introduction}
Fulton asked how many solutions to a problem of enumerative
geometry can be real, when that problem is one of counting geometric figures
of some kind having specified position with respect to some general fixed
figures~\cite{Fu_84}.
For the problem of plane conics tangent to five general conics, the
(surprising) answer is that all 3264 may be real~\cite{RTV}. 
Similarly, given any problem of enumerating $p$-planes incident on
some general fixed subspaces, there are real fixed
subspaces such that each of the (finitely many) incident  $p$-planes are
real~\cite{So99}. 
We show that the problem of enumerating parameterized rational curves in a
Grassmannian satisfying simple (codimension 1) conditions may have all of
its solutions be real.

This problem of enumerating rational curves on a Grassmannian arose in at
least two distinct areas of mathematics.
The number of such curves was predicted by the formula of Vafa and
Intriligator~\cite{Vafa,Intriligator} from mathematical physics.
It is also the number of complex dynamic compensators which stabilize a
particular linear system, and the enumeration was solved in this
context~\cite{RRW98,RRW96}.  
The question of real solutions also arises in systems theory~\cite{Byrnes}.
Our proof, while exploiting techniques from systems theory, has no direct
implications for the problem of real dynamic output compensation.

\section{Statement of results}

We work with complex algebraic varieties and ask when {\it a
priori}$\,$ complex solutions to an enumerative problem are real.
Fix integers $m,p>1$ and $q\geq 0$.
Set $n:=m+p$.
Let ${\bf G}$ be the Grassmannian of $p$-planes in 
${\mathbb C}^n$.
The space ${\mathcal M}_q$ of maps 
$M:{\mathbb P}^1\rightarrow{\bf G}$ of degree $q$ has
dimension $N:=pm+qn$~\cite{Clark,Stromme}.
If $L$ is an $m$-plane and $s\in {\mathbb P}^1$, then the collection of
all maps $M$ satisfying $M(s)\cap L\neq \{0\}$ is an
irreducible subvariety of codimension 1.
We study the following enumerative problem:
\begin{equation}
  \label{eq:enumerative}
  \begin{minipage}{5.2in}
    Given general points $s_1,\ldots,s_N$ in ${\mathbb P}^1$ and general 
    $m$-planes $L_1,\ldots,L_N$ in ${\mathbb C}^n$, 
    how many maps $M\in {\mathcal M}_q$
    satisfy $M(s_i)\cap L_i\neq\{0\}$ for $i=1,\ldots,N$?
  \end{minipage}
\end{equation}

Rosenthal~\cite{Rosen94} interpreted the solutions as a linear section of a
projective embedding of ${\mathcal M}_q$, and 
Ravi, Rosenthal, and Wang~\cite{RRW98,RRW96} show that 
the degree of its closure ${\mathcal K}_q$ in this embedding is
\begin{equation}\label{dqmp}
  \delta\quad:=\quad(-1)^{q(m+1)}N! 
   \sum_{\nu_1+\cdots+\nu_p=q}
  \frac{\prod_{i<j}(j-i+n(\nu_j-\nu_i))}
  {\prod_{j=1}^p(m+j+n\nu_j-1)!}\ .
\end{equation}
Thus, if there are finitely many solutions, then their number (counted with 
multiplicity) is at most $\delta$.
The difference between $\delta$ and the number of solutions counts points
common to both the linear section and the boundary 
${\mathcal K}_q-{\mathcal M}_q$ of ${\mathcal K}_q$.
Since ${\bf G}$ is a homogeneous space, an application of
Kleiman's Theorem~\cite{Kleiman} shows there are finitely many solutions
and no multiplicities.
Bertram~\cite{Bertram} uses explicit methods (a moving lemma) to to show
there are finitely many solutions and also no points in the boundary of 
${\mathcal Q}_q$, and hence none in the boundary of ${\mathcal K}_q$.
He also computes the small quantum cohomology ring of ${\bf G}$, which gives
algorithms for computing $\delta$ and other intersection 
numbers involving rational curves on a Grassmannian.

When the $s_i$ and $L_i$ are real, not all of these solutions are defined
over the real numbers.
We show there are real $s_i$ and $L_i$ for which each of the
$\delta$ maps are real.

\begin{thm}\label{thm:real}
There exist real $m$-planes $L_1,\ldots,L_N$ in ${\mathbb R}^n$ and
points  $s_1,\ldots,s_N\in{\mathbb P}^1_{\mathbb R}$ so that 
there are exactly $\delta$ maps
$M:{\mathbb P}^1\rightarrow{\bf G}$ of degree $q$
which satisfy $M(s_i)\cap L_i \neq \{0\}$ for each $i=1,\ldots,N$,
and each of these are real.
\end{thm}

Our proof is elementary in that it argues from the equations for the
locus of maps $M$ which satisfy $M(s)\cap L\neq \{0\}$.
A consequence is that we obtain fairly explicit choices of $s_i$ and $L_i$
which give only real maps, which we discuss in Section 4.
Also, our proof uses neither Kleiman's Theorem nor Bertram's moving
lemma, and thus it provides a new and elementary proof that there are
$\delta$ solutions to the enumerative problem~(\ref{eq:enumerative}).

\section{The quantum Grassmannian}

The space ${\mathcal M}_q$ of maps 
${\mathbb P}^1\rightarrow{\bf G}$ of degree $q$ is a smooth
quasi-projective algebraic variety.
A smooth compactification is provided by a quot scheme 
${\mathcal Q}_q$~\cite{Stromme}.
By definition, there is a universal exact sequence
$$
  0\ \rightarrow\ {\mathcal S}\ \rightarrow\
  {\mathbb C}^n\otimes{\mathcal O}\ \rightarrow\ 
  {\mathcal T}\ \rightarrow 0
$$
of sheaves on ${\mathbb P}^1\times{\mathcal Q}_q$ where 
${\mathcal S}$ is a vector bundle of degree $-q$ and rank $p$.
Twisting the determinant of ${\mathcal S}$ by 
${\mathcal O}_{{\mathbb P}^1}(q)$ and pushing forward to 
${\mathcal Q}_q$ induces a Pl\"ucker map
$$
  {\mathcal Q}_q\ \rightarrow\ 
  {\textstyle {\mathbb P}\left( \bigwedge^p{\mathbb C}^n\otimes 
      H^0({\mathcal O}_{{\mathbb P}^1}(q))^*\right)}
$$
which is the analog of the Pl\"ucker embedding of ${\bf G}$.
The Pl\"ucker map is an embedding of ${\mathcal M}_q$, and so its image 
${\mathcal K}_q$ provides a different compactification of 
${\mathcal M}_q$.
We call ${\mathcal K}_q$ the quantum Grassmannian.
(In~\cite{BDW}, this space is called the Uhlenbeck compactification).
Our proof of Theorem~\ref{thm:real} exploits some of its structures that
were elucidated in work in systems theory.

The Pl\"ucker map fails to be injective on the boundary 
${\mathcal Q}_q-{\mathcal M}_q$ of ${\mathcal Q}_q$.
Indeed, Bertram~\cite{Bertram} constructs a ${\mathbb P}^{p-1}$ bundle over 
${\mathbb P}^1\times{\mathcal Q}_{q-1}$ that maps onto the boundary, with its
restriction over ${\mathbb P}^1\times{\mathcal M}_{q-1}$ an embedding.
On this projective bundle, the Pl\"ucker map factors through the base 
${\mathbb P}^1\times{\mathcal Q}_{q-1}$ and the image of a point in the base
is $s\cdot S$, where $s$ is the section of 
${\mathcal O}_{{\mathbb P}^1}(1)$ vanishing at $s\in{\mathbb P}^1$ and 
$S$ is the image of a point in ${\mathcal Q}_{q-1}$ under its Pl\"ucker map.
This identifies the image of the exceptional locus of the Pl\"ucker map 
with the image of ${\mathbb P}^1\times{\mathcal K}_{q-1}$ in 
${\mathcal K}_q$ under a map $\pi$ (given below).

More concretely, a
point in ${\mathcal Q}_q$ may be  (non-uniquely) represented by
a $p\times n$-matrix $M$ of forms in $s,t$, with
homogeneous rows and whose maximal minors have degree $q$~\cite{RR94}.
The image of such a point under the Pl\"ucker map is the collection of
maximal minors of $M$.
The maps in ${\mathcal M}_q$ are represented by matrices whose maximal
minors have no common factors:
Given such a matrix $M$, the association
$$
{\mathbb P}^1\ni(s,t)\ \longmapsto\ \mbox{row space }M(s,t)
$$
defines a map of degree $q$.

The collection $\binom{[n]}{p}$ of $p$-subsets of $\{1,\ldots,n\}$
index the maximal minors of $M$.
For $\alpha\in\binom{[n]}{p}$ and $0\leq a\leq q$, the coefficients
$z_{\alpha^{(a)}}$ of $s^at^{q-a}$ in  
the $\alpha$th maximal minor of $M$ provide Pl\"ucker coordinates
for maps in ${\mathcal M}_q$, and for the space 
${\mathbb P}\left(\bigwedge^p{\mathbb C}^n\otimes 
      H^0({\mathcal O}_{{\mathbb P}^1}(q))^*\right)$.
Let ${\mathcal C}_q:=
     \{\alpha^{(a)}\mid \alpha\in\binom{[n]}{p}, 0\leq a\leq q\}$ be the
indices of these Pl\"ucker coordinates.
Then the image of the exceptional locus in 
${\mathcal K}_q$ is the image of the (birational) map
$\pi: {\mathbb P}^1\times{\mathcal K}_{q-1}\rightarrow{\mathcal K}_q$
defined by 
\begin{equation}\label{eq:pi}
  \pi\ :\  \left([A,B],(x_{\beta^{(b)}}\mid 
   \beta^{(b)}\in{\mathcal C}_{q-1})\right)
   \ \longmapsto\ 
  (Ax_{\alpha^{(a)}}-Bx_{\alpha^{(a-1)}}\mid
      \alpha^{(a)}\in{\mathcal C}_q)\,.
\end{equation}

\rule{0pt}{15pt}
The relevance of the quantum Grassmannian ${\mathcal K}_q$ to the enumerative
problem~(\ref{eq:enumerative}) is seen by considering the condition 
for a map $M\in{\mathcal M}_q$ to satisfy $M(s,t)\cap L\neq\{0\}$
where $L$ is an $m$-plane in ${\mathbb C}^n$ and $(s,t)\in{\mathbb P}^1$.
If we represent $L$ as the row space of a $m\times n$ matrix, also written
$L$, then this condition is
$$
  0\ =\ 
  \det\left[\begin{array}{c}L\\M(s,t)\end{array}\right]\ =\ 
  \sum_{\alpha\in\binom{[n]}{p}} f_\alpha(s,t)\, l_\alpha\,,
$$
the second expression given by Laplace expansion of the
determinant along the rows of $M$.
Here, $l_\alpha$ is the appropriately signed maximal minor of $L$.
If we expand the forms $f_\alpha(s,t)$ in this last expression, we obtain
$$
  \sum_{\alpha^{(a)}\in{\mathcal C}_q} 
    z_{\alpha^{(a)}} s^a t^{q-a}l_\alpha\ =\ 0\,,
$$
a linear equation in the Pl\"ucker coordinates of $M$.
Thus the solutions $M\in{\mathcal M}_q$ to the enumerative
problem~(\ref{eq:enumerative}) are a linear section of ${\mathcal M}_q$ in
its Pl\"ucker embedding, and so the degree $\delta$ of ${\mathcal K}_q$
provides an upper bound on the number of solutions.

\begin{figure}[htb]
 $$\epsfxsize=1.8in \epsfbox{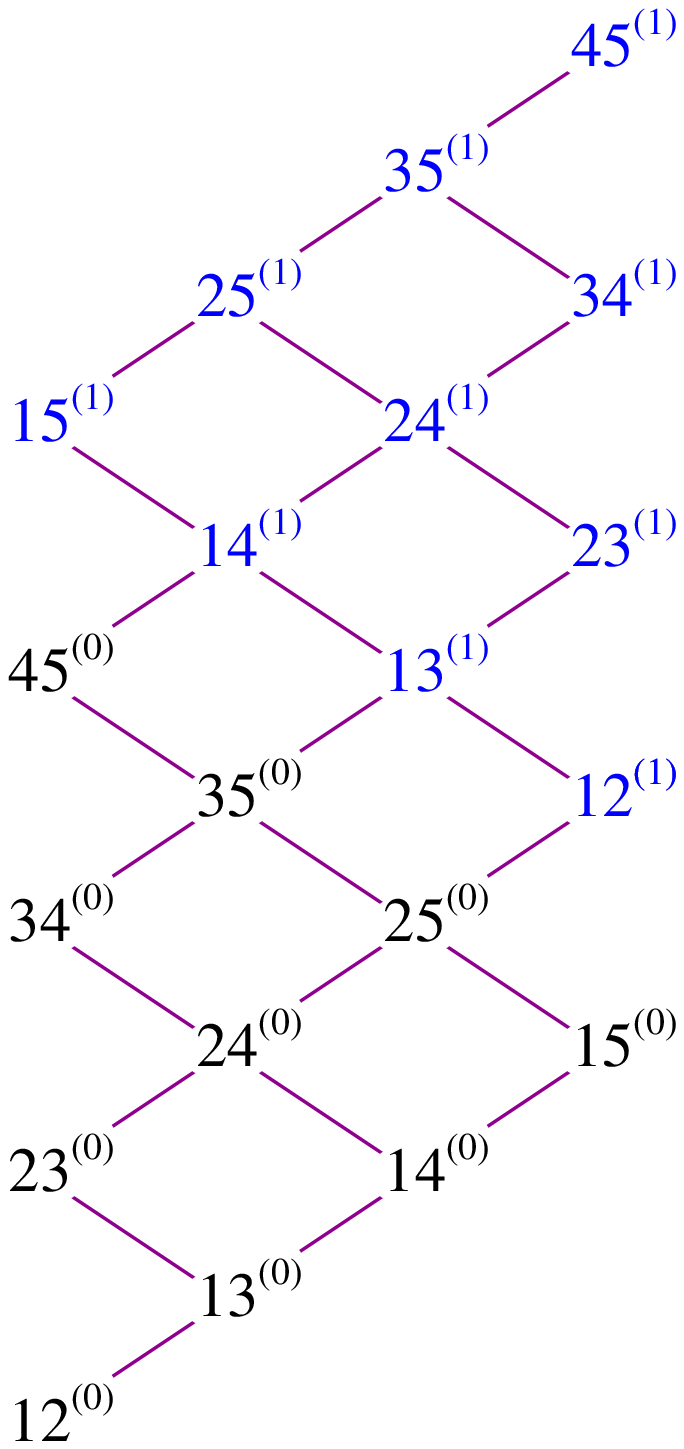}$$
 \caption{${\mathcal C}_1$\label{fig:one}.}
\end{figure}

The set ${\mathcal C}_q$ of Pl\"ucker coordinates has a natural partial
order 
$$
  \alpha^{(a)}\ \geq\ \beta^{(b)} \quad\Longleftrightarrow\quad
  \begin{array}{c}
   a\geq b, \mbox{ and if } a-b<p, \mbox{ then }\\
   \alpha_{a-b+1}\geq \beta_1, 
   \ldots,\alpha_p\geq \beta_{p+1-b+a}
  \end{array}\ .
$$
The poset ${\mathcal C}_q$ is graded with the rank, $|\alpha^{(a)}|$, of
$\alpha^{(a)}$ equal to $an + \sum_i \alpha_i-i$.
Figure~\ref{fig:one} shows ${\mathcal C}_1$ when $p=2$ and $m=3$.
Given $\alpha^{(a)}\in {\mathcal C}_q$, define the {\it quantum
Schubert variety}
$$
  Z_{\alpha^{(a)}}\quad :=\quad \{z=(z_{\beta^{(b)}})\in {\mathcal K}_q\mid
  z_{\beta^{(b)}}=0\mbox{ if }\ \beta^{(b)}\not\leq \alpha^{(a)} \}\ .
$$

Let ${\mathcal H}_{\alpha^{(a)}}$ be the hyperplane 
defined by $z_{\alpha^{(a)}}=0$.
The main technical result we use is the following.

\begin{prop}[\cite{RRW96,RRW98}]\label{prop:RRW}
Let $\alpha^{(a)}\in{\mathcal C}_q$.
Then
\begin{enumerate}
\item[(i)]
   $Z_{\alpha^{(a)}}$ is an irreducible subvariety of $\,{\mathcal K}_q$ of
   dimension $|\alpha^{(a)}|$.
\item[(ii)]
   The intersection of $\,Z_{\alpha^{(a)}}$ and 
   ${\mathcal H}_{\alpha^{(a)}}$ is
   generically transverse, and 
   $$
     Z_{\alpha^{(a)}}\cap {\mathcal H}_{\alpha^{(a)}}\ =\ 
     \bigcup_{\beta^{(b)}\lessdot\alpha^{(a)}} Z_{\beta^{(b)}}\,.
   $$
\end{enumerate}
\end{prop}

Another proof of (ii) is given in~\cite{SS_SAGBI}, which shows (ii) is an
ideal-theoretic equality.
From (ii) and B\'ezout's theorem, we obtain the following recursive formula
for the degree of $Z_{\alpha^{(a)}}$
$$
\deg Z_{\alpha^{(a)}}\  =\ 
       \sum_{\beta^{(b)}\lessdot\alpha^{(a)}} \deg Z_{\beta^{(b)}}\,.
$$
Since the minimal quantum Schubert variety is a point, we deduce the main
result of~\cite{RRW98}:

\begin{cor}
The degree $\delta$ of $\,{\mathcal K}_q$  is the number of
maximal chains in the poset ${\mathcal C}_q$.
\end{cor}

Closed formulas are given for $\delta$ in~\cite{RRW96,RRW98}, 
the source of the formula~(\ref{dqmp}), as well as the
number $\deg Z_{\alpha^{(a)}}$ of maximal chains below $\alpha^{(a)}$.

\section{Proof of Theorem~\ref{thm:real}}

Let $L(s,t)$ be the $m$-plane osculating the parameterized rational normal
curve 
$$
\gamma\ :\ (s,t)\in{\mathbb P}^1\  \longmapsto\ 
 (s^{n-1},\,ts^{n-2},\,\ldots,\,t^{n-2}s,\,t^{n-1})\in{\mathbb P}^{n-1}
$$
at the point $\gamma(s,t)$.
Then $L(s,t)$ is the row space of the $m\times n$ matrix of forms with rows
$\gamma(s,t),\gamma'(s,t),\ldots,\gamma^{(m-1)}(s,t)$, the derivative taken 
with respect to the parameter $t$.
Write $L(s,t)$ for this matrix.
For $\alpha\in\binom{[n]}{p}$, the maximal minor of $L(s,t)$
complementary to $\alpha$ is 
$(-1)^{|\alpha|}s^{\binom{m}{2}}l_\alpha s^{|\alpha|}t^{mp-|\alpha|}$,
where $|\alpha|:=\sum_i \alpha_i-i$ and $(-1)^{|\alpha|}l_\alpha$ is the
corresponding maximal minor of $L(1,1)$.
Let ${\mathcal H}(s,t)$ be the pencil of hyperplanes 
given by the linear form
$$
  \Lambda(s,t)\ :=\ \sum_{\alpha^{(a)}\in{\mathcal C}_q} 
    z_{\alpha^{(a)}}
     l_\alpha s^{|\alpha^{(a)}|}t^{N-|\alpha^{(a)}|}\,.
$$
Let $M$ be a matrix representing a curve in ${\mathcal M}_q$.
Then
$$
  \det\left[\begin{array}{c} 
        L(s,t)\\M(s^n,t^n)\rule{0pt}{16pt}\end{array}\right]\ =\ 
  s^{\binom{m}{2}} \sum_{\alpha^{(a)}\in{\mathcal C}_q} z_{\alpha^{(a)}}
  s^{an}t^{(q-a)n} l_\alpha s^{|\alpha|}t^{mp-|\alpha|}\ =\    
  s^{\binom{m}{2}}\Lambda(s,t)\,.
$$
Thus ${\mathcal M}_q\cap {\mathcal H}(s,t)$ consists of all maps 
$M:{\mathbb P}^1\rightarrow {\bf G}$ of degree $q$ which satisfy
$M(s^n,t^n)\cap L(s,t)\neq \{0\}$.

Theorem~1 is a consequence of the following two theorems.

\begin{thm}\label{thm:transverse}
There exist positive real numbers $t_1,\ldots,t_N$ such that for any
$\alpha^{(a)}\in{\mathcal C}_q$, the intersection
$$
 Z_{\alpha^{(a)}}\cap{\mathcal H}(1,t_1)\cap\cdots\cap
{\mathcal H}(1,t_{|\alpha^{(a)}|})
$$ 
is transverse with all points of intersection real.
\end{thm}

\begin{thm}\label{thm:proper}
If $t_1,\ldots,t_k\in {\mathbb C}$ are distinct, 
then for any $\alpha^{(a)}\in{\mathcal C}_q$,
the intersection  
\begin{equation}\label{eq:proper}
  Z_{\alpha^{(a)}}\cap{\mathcal H}(1,t_1)\cap\cdots\cap{\mathcal H}(1,t_k)
\end{equation}
is proper in that it has dimension $|\alpha^{(a)}|-k$.
\end{thm}

\noindent{\bf Proof of Theorem~\ref{thm:real}. }
By Theorem~\ref{thm:transverse}, there exist positive real numbers
$t_1,\ldots,t_N$ (necessarily distinct) so that the intersection 
\begin{equation}\label{eq:int}
  {\mathcal K}_q \cap 
     {\mathcal H}(1,t_1)\cap\cdots\cap{\mathcal H}(1,t_N)
\end{equation}
is transverse and consists of exactly $\delta$ real points.
We show all these points lie in ${\mathcal M}_q$, and thus are maps
$M:{\mathbb P}^1\rightarrow {\bf G}$ of degree $q$ satisfying
$M(1,t_i^n)\cap L(1,t_i)\neq\{0\}$ for $i=1,\ldots,N$, which proves
Theorem~\ref{thm:real}. 

Recall the map 
$\pi:{\mathbb P}^1\times{\mathcal K}_{q-1}\rightarrow
{\mathcal K}_q$~(\ref{eq:pi}) whose
image is the complement of ${\mathcal M}_q$ in ${\mathcal K}_q$.
Then 
\begin{eqnarray*} 
  \pi^*{\mathcal H}(s,t)&=&
   \sum_{\alpha^{(a)}\in{\mathcal C}_q}
    (Ax_{\alpha^{(a)}}-Bx_{\alpha^{(a-1)}}) 
       l_\alpha s^{|\alpha^{(a)}|}t^{N-|\alpha^{(a)}|}\\
  &=&
    (At^n-Bs^n)\sum_{\alpha^{(a)}\in{\mathcal C}_{q-1}}
    x_{\alpha^{(a)}}l_\alpha s^{|\alpha^{(a)}|}t^{N-n-|\alpha^{(a)}|}\,.
\end{eqnarray*}
Hence, if ${\mathcal H}'(s,t)$ is the pencil of hyperplanes in the Pl\"ucker
space of ${\mathcal K}_{q-1}$ defining the locus of $M\in{\mathcal M}_{q-1}$
satisfying $M(s^n,t^n)\cap L(s,t)\neq \{0\}$, then 
$$
  \pi^*{\mathcal H}(s,t)\ =\ 
  (At^n-Bs^n){\mathcal H}'(s,t)\,.
$$

Thus any point in~(\ref{eq:int}) not in ${\mathcal M}_q$ is the
image of a point 
$([A,B],M)$ in ${\mathbb P}^1\times{\mathcal K}_{q-1}$
satisfying 
$\pi^*{\mathcal H}(1,t_i)=(At_i^n-B){\mathcal H}'(1,t_i)$ 
for each $i=1,\ldots,N$.
As the $t_i$ are positive and distinct, 
such a point can only satisfy $At_i^n-B=0$ for one $i$.
Thus $M\in{\mathcal K}_{q-1}$ lies in 
at least $N-1$ of the hyperplanes ${\mathcal H}'(1,t_i)$.
Since $N-1$ exceeds the dimension $N-n$ of ${\mathcal K}_{q-1}$, there are
no such 
points $M\in{\mathcal K}_{q-1}$, by Theorem~\ref{thm:proper} for maps of
degree $q-1$.
\qed\bigskip

\noindent{\bf Proof of Theorem~\ref{thm:proper}. }
For any $t_1,\ldots,t_k$, the intersection~(\ref{eq:proper}) has dimension 
at least $|\alpha^{(a)}|-k$.
We show it has at most this dimension, if $t_1,\ldots,t_k$ are distinct.

Suppose $k=|\alpha^{(a)}|+1$ and let $z\in  Z_{\alpha^{(a)}}$.
Then $z_{\beta^{(b)}}=0$ if $\beta^{(b)}\not\leq\alpha^{(a)}$
and so the form $\Lambda(1,t)(z)$ defining ${\mathcal H}(1,t)$ is divisible
by $t^{N-|\alpha^{(a)}|}$ with quotient
$$
\sum_{\beta^{(b)}\leq\alpha^{(a)}} z_{\beta^{(b)}}l_\beta\,
t^{|\alpha^{(a)}|-|\beta^{(b)}|}\,.
$$
This is a non-zero polynomial in $t$ of degree at most $|\alpha^{(a)}|$
and thus it vanishes for at most $|\alpha^{(a)}|$ distinct $t$.
It follows that~(\ref{eq:proper}) is empty for $k>|\alpha^{(a)}|$.

If $k\leq |\alpha^{(a)}|$ and $t_1,\ldots,t_k$ are distinct,
but~(\ref{eq:proper}) has dimension exceeding  
$|\alpha^{(a)}|-k$, then completing $t_1,\ldots,t_k$ to a set of distinct 
numbers $t_1,\ldots,t_{|\alpha^{(a)}|+1}$ would give a non-empty
intersection in~(\ref{eq:proper}), a contradiction.
\qed\bigskip

\noindent{\bf Proof of Theorem~\ref{thm:transverse}. }
We construct the sequence $t_i$ inductively.
If we let $\alpha=1<2<\cdots<p-1<p+1$, then 
$ Z_{\alpha^{(0)}}$ is a line.
Indeed, it is isomorphic to the set of $p$-planes containing a fixed
$(p-1)$-plane and lying in a fixed $(p+1)$-plane.
By Theorem~\ref{thm:proper}, $Z_{\alpha^{(0)}}\cap{\mathcal H}(1,t)$ is
then a single, necessarily real, point, for any real number $t$.
Let $t_1$ be any positive real number.

Suppose we have positive real numbers $t_1,\ldots,t_k$ with the property
that for any $\beta^{(b)}$ with $|\beta^{(b)}|\leq k$, 
$$
 Z_{\beta^{(b)}}\cap{\mathcal H}(1,t_1)\cap\cdots\cap
{\mathcal H}(1,t_{|\beta^{(b)}|})
$$ 
is transverse with all points of intersection real.

Let $\alpha^{(a)}$ be an index with $|\alpha^{(a)}|=k+1$ and consider the
1-parameter family ${\mathcal Z}(t)$ of schemes defined for $t\neq 0$ by 
$ Z_{\alpha^{(a)}}\cap{\mathcal H}(1,t)$.
For $t\neq 0$, if we restrict the form $\Lambda(1,t)$ to 
$z\in Z_{\alpha^{(a)}}$, then, after dividing out $t^{N-|\alpha^{(a)}|}$, 
we obtain
$$
z_{\alpha^{(a)}} + \sum_{\beta^{(b)}<\alpha^{(a)}} z_{\beta^{(b)}}l_\beta\, 
t^{|\alpha^{(a)}|-|\beta^{(b)}|}\,.
$$
Thus ${\mathcal Z}(0)$ is
$$
    Z_{\alpha^{(a)}}\cap {\mathcal H}_{\alpha^{(a)}}\ =\ 
   \bigcup_{\beta^{(b)}\lessdot\alpha^{(a)}}  Z_{(\beta,d)}\,,
$$
by Proposition~\ref{prop:RRW} (ii).

\noindent{\bf Claim:}
The cycle
$$
   {\mathcal Z}(0)\cap{\mathcal H}(1,t_1)\cap\cdots\cap{\mathcal H}(1,t_k)
$$
is free of multiplicities.

If not, then there are two components $ Z_{\beta^{(b)}}$
and $ Z_{\gamma^{(c)}}$ of ${\mathcal Z}(0)$ such that 
$$
  Z_{\beta^{(b)}} \cap Z_{\gamma^{(c)}}
  \cap {\mathcal H}(1,t_1)\cap \cdots\cap {\mathcal H}(1,t_k)
$$
is nonempty.
But this contradicts Theorem~\ref{thm:proper}, as 
$Z_{\beta^{(b)}} \cap Z_{\gamma^{(c)}}=Z_{\delta^{(d)}}$, where
$\delta^{(d)}$ is the greatest lower bound of $\beta^{(b)}$ and
$\gamma^{(c)}$ in ${\mathcal C}_q$, and so  
$\dim Z_{\delta^{(d)}}<\dim Z_{\beta^{(b)}}=k$.

From the claim, there is an $\epsilon_{\alpha^{(a)}}>0$ such that if
$0\leq t\leq \epsilon_{\alpha^{(a)}}$, then 
$$
  {\mathcal Z}(t)\cap{\mathcal H}(1,t_1)\cap\cdots\cap{\mathcal H}(1,t_k)
$$
is transverse with all points of intersection real.
Set
$$
  t_{k+1}\ :=\ \min\{\epsilon_{\alpha^{(a)}} : |\alpha^{(a)}|=k+1\}\,.
   \qquad\qed
$$

\section{Further Remarks}
From our proof of Theorem~\ref{thm:transverse}, we obtain a rather precise
choice of $s_i$ and $L_i$ in the enumerative problem which give only real
maps.
By $\forall t_1\gg t_2\gg\cdots\gg t_N>0$, we mean
$$
\forall t_1>0\ \ \exists \epsilon_2>0\ \ \forall \epsilon_2>t_2>0\ \ 
\cdots\ \ \exists \epsilon_N>0\ \ \forall \epsilon_N>t_N>0\, .
$$

\begin{cor}
 $\forall t_1\gg t_2\gg\cdots\gg t_N>0$, each of the $\delta$ maps
$M:{\mathbb P}^1\rightarrow{\bf G}$ of degree $q$ which satisfy
$M(1,t_i)\cap L(1,t_i^{1/n})\neq \{0\}$
for $i=1,\ldots,N$ are real.
\end{cor}

When $q=0$, there is substantial evidence~\cite{Sottile_shapiro}
that this choice of $t_1,\ldots,t_N$ is too restrictive.
B.~Shapiro and M.~Shapiro have the following conjecture:
\medskip

\noindent{\bf Conjecture. }
{\it 
Suppose $q=0$.
Then for generic real numbers $\,t_1,\ldots,t_{mp}$ all of the finitely many  
$p$-planes $H$ which satisfy $H\cap L(1,t_i)\neq \{0\}$ are real.
}\medskip

In contrast, when $q>0$, the restriction 
$\forall t_1\gg t_2\gg\cdots\gg t_N>0$ is necessary.
We observe this in the case when $q=1$, $p=m=2$, so $N=8$ and $\delta=8$.
That is, for parameterized curves of degree 1 in the Grassmannian of
2-planes in ${\mathbb C}^4$.
Here, the choice of $t_i=i$ in~(\ref{eq:int}) gives no real maps, while
the choice  $t_i=i^6$ gives 8 real maps.

We briefly describe that calculation.
There are 12 Pl\"ucker coordinates $z_{ij^{(a)}}$ for $1\leq i<j\leq 4$ and
$a=0,1$. 
If we let $f_{ij}:=tz_{ij^{(0)}}+sz_{ij^{(1)}}$, then 
$$
f_{14}f_{23} - f_{13}f_{24} + f_{12}f_{34}\ =\ 0\,,
$$
as $f_{ij}(s,t)\in{\bf G}$ for all $s,t$.
The coefficients of $t^2$, $st$, and $s^2$ in this expression give three
quadratic relations among the $z_{ij^{(a)}}$:
$$ 
\begin{array}{c}
 z_{14^{(0)}}z_{23^{(0)}}\ -\ 
   z_{13^{(0)}}z_{24^{(0)}}\ +\ z_{12^{(0)}}z_{34^{(0)}}\,,\\
 z_{12^{(1)}}z_{34^{(0)}}\ -\ 
   z_{13^{(1)}}z_{24^{(0)}}\ +\  z_{14^{(1)}}z_{23^{(0)}} \rule{0pt}{15pt}
   + z_{23^{(1)}}z_{14^{(0)}}\ -\ 
   z_{24^{(1)}}z_{13^{(0)}}\ +\ z_{34^{(1)}}z_{12^{(0)}}\,,\\
 z_{14^{(1)}}z_{23^{(1)}}\ -\            \rule{0pt}{15pt}
   z_{13^{(1)}}z_{24^{(1)}}\ +\ z_{12^{(1)}}z_{34^{(1)}}\,,
\end{array}
$$ 
and these constitute a Gr\"obner basis for the
homogeneous ideal of ${\mathcal K}_1$\cite{SS_SAGBI}.

Here, the form $\Lambda$ is
$$
\begin{array}{c}
{\displaystyle 
\ \ 
t^8z_{12^{(0)}}-2t^7z_{13^{(0)}}+t^6z_{14^{(0)}}+3t^6z_{23^{(0)}}
-2t^5z_{24^{(0)}}+t^4z_{34^{(0)}}}\\
+\,t^4z_{12^{(1)}}-2t^3z_{13^{(1)}}+t^2z_{14^{(1)}}+3t^2z_{23^{(1)}}
-2t\ \, z_{24^{(1)}}+\ \ z_{34^{(1)}}\,.
\end{array}
$$
We set $z_{34^{(1)}}=1$ and work in local coordinates.
Then the ideal generated by the 3 quadratic equations 
and 8 linear relations $\Lambda(t_i)$ for $i=1,\ldots,8$
defines the 8 solutions to~(\ref{eq:int}).
We used  Maple V.5 to generate these equations and then  compute 
a univariate polynomial in the ideal, which had degree 8.
This polynomial had no real solutions when $t_i=i$, but all 8 were
real when $t_i=i^6$.
(Elimination theory guarantees that the number of real solutions equals the
number of real roots of the eliminant.)
\medskip

We describe how the enumerative problem~(\ref{eq:enumerative}) arises in
systems theory (see also~\cite{Byrnes}).
A physical system (eg.~a mechanical linkage) with $m$ inputs and $p$
measured outputs whose evolution is governed by a system of linear
differential equations is modeled by a $m\times n$-matrix $L(s)$ of
real univariate polynomials.
The largest degree of a maximal minor of this matrix is the 
MacMillan degree, $r$, of the evolution equation.
Consider now controlling this linear system by output feedback with a
dynamic compensator.
That is, a $p$-input, $m$-output linear system $M$ is used to couple the $m$ 
inputs of the system $L$ to its $p$ outputs.
The resulting closed system has characteristic polynomial
$$
\varphi(s)\ :=\ \left[\begin{array}{c}L(s)\\M(s)\end{array}\right]\,,
$$
and the roots of $\varphi$ are the natural frequencies or {\it poles} of the
closed system.
The dynamic pole assignment problem asks, given a system $L(s)$ and a
desired characteristic polynomial $\varphi$, can one find a (real)
compensator $M(s)$ of MacMillan degree $q$ so that the resulting closed
system has characteristic polynomial $\varphi$?
That is, if $s_1,\ldots,s_{r+p}$ are the roots of $\varphi$, which
$M\in{\mathcal M}_q$ satisfy 
$$
\det\left[\begin{array}{c}L(s_i)\\M(s_i)\end{array}\right]\ =\ 0,
\qquad\mbox{for }i=1,2,\ldots,r+p\,?
$$

In the critical case when $r+q=mp+qn$, this is an instance of the
enumerative problem~(\ref{eq:enumerative}).
When the degree $\delta$ is odd, then for a real system $L$ and a real
characteristic polynomial $\varphi$, there will be at least one real dynamic
compensator. 
Part of the motivation for~\cite{RRW96} was to obtain a formula for $\delta$ 
from which its parity could be deduced for different values of
$q,m$, and $p$. 

From this description, we see that the choice of planes $L_i$ that arise
in the dynamic pole placement problem are 
$N=mp+qn$ points on a rational curve of degree $mp+(n-1)q$ in the
Grassmannian of $m$-planes in ${\mathbb C}^n$.
In contrast, the planes of Theorem~\ref{thm:transverse} (and hence of
Theorem~\ref{thm:real}) arise as $N$ points on a 
rational curve of degree $mp$.
Only when $q=0$ (the case of static compensators) is there any overlap.
While our proof of Theorem~\ref{thm:real} owes much to systems theory, it
has no direct implications for the problem of real dynamic output
compensation. 
\medskip

Our method of proof of Theorem~\ref{thm:real} (like that in~\cite{So99}) was
inspired by the numerical 
Pieri homotopy algorithm of~\cite{HSS} for computing the solutions
to~(\ref{eq:enumerative})  when $q=0$.
Likewise, the explicit degenerations of intersections of the 
${\mathcal H}(s,t)$ that we used, and more generally
Proposition~\ref{prop:RRW} (ii), can be used to construct an optimal
numerical homotopy algorithm for finding the solutions
to~(\ref{eq:enumerative}).
This is in exactly the same manner as the explicit degenerations of
intersections of special Schubert varieties of~\cite{So97d} were used to
construct the Pieri homotopy algorithm of~\cite{HSS}.
\medskip

We close with one open problem concerning the enumeration of
rational curves on a Grassmannian.
For a point $s\in{\mathbb P}^1$ and any Schubert variety $\Omega$ of 
${\bf G}$, consider the quantum Schubert variety $\Omega(s)$ of curves
$M\in{\mathcal M}_q$ satisfying $M(s)\in\Omega$.
The quantum Schubert calculus gives algorithms to compute the number of
curves $M\in{\mathcal M}_q$ which lie in the intersection of an appropriate
number of these $\Omega(s)$, and we ask when it is possible to have all
solutions real.
A modification of the proof of Theorem~\ref{thm:transverse} shows that
this is the case when all except possibly 2 are hypersurface Schubert
varieties.
In every case we have been able to compute, all solutions may be real.

\end{document}